\documentstyle[12pt,german]{article}
\begin{document}
\includegraphics{zibheader.eps}
\vspace*{6.5cm}
\begin{center}
{\Large Wolfram Koepf}\vspace*{2cm}\\
\LARGE{{\bf Algebraische Darstellung}}\\[2mm]
\LARGE{{\bf transzendenter Funktionen}}
\end{center}
\vfill
\hrule
\vspace*{3mm}
Preprint SC 94--24 (Oktober 1994)

\thispagestyle{empty}
\setcounter{page}{0}
\eject
\hoffset -1cm
\voffset -2cm
\footskip=1cm
\parindent=0pt

\newcommand\tm{$^{\rm tm}$}
\def\CA{Com\-pu\-ter\-al\-ge\-bra}
\def\ca{Com\-pu\-ter\-al\-ge\-bra}
\def\cas{Com\-pu\-ter\-al\-ge\-bra-Sy\-stem}
\def\bibtex{{\sc Bib}\TeX}

\newcommand\voidbox[1]{#1}
\newcommand\elib{eLib}
\newcommand\cais{{\sc CAIS}}
\newcommand\infocais{Genauere Informationen finden Sie im \cais.}

\newcommand{\axiom}{{\sc AXIOM}}
\newcommand{\derive}{{\sc Derive}}
\newcommand{\macsyma}{{\sc Macsyma}}
\newcommand{\maple}{{\sc Maple}}
\newcommand{\mathematica}{{\sc Mathematica}}
\newcommand{\mathem}{{\sc Mathematica}}
\newcommand{\reduce}{{\sc Reduce}}
\newcommand{\aldessac}{{\sc Aldes/SAC-2}}
\newcommand{\sac}{{\sc SAC-2}}
\newcommand{\aldes}{{\sc ALDES}}
\newcommand{\cayley}{{\sc Cayley}}
\newcommand{\magma}{{\sc MAGMA}}
\newcommand{\cocoa}{{\sf CoCoA}}
\newcommand{\felix}{{\sc FELIX}}
\newcommand{\form}{{\sc FORM}}
\newcommand{\gap}{{\sf GAP}}
\newcommand{\kant}{{\sc KANT-V2}}
\newcommand{\mas}{{MAS}}
\newcommand{\masyca}{{MASYCA}}
\newcommand{\mathcad}{{\sc Mathcad} }
\newcommand{\macaulay}{{\sc Macaulay}}
\newcommand{\macauley}{{\sc Macaulay}}
\newcommand{\mupad}{{\rm MuPAD}}
\newcommand{\chevie}{{\sf CHEVIE}}
\newcommand{\qpic}{{\sc Quotpic}}
\newcommand{\pari}{{PARI}}
\newcommand{\reve}{{\sc Reve}}
\newcommand{\orme}{{\sc ORME}}
\newcommand{\simath}{{\sc Simath}}
\newcommand{\parsac}{{\sc Parsac}}
\newcommand{\symmetrica}{{\sc Symmetrica}}
\newcommand{\reptiles}{{\sc RepTiles}}
\newcommand{\symmpad}{{\sc SymmPAD}}
\newcommand{\singular}{{SINGULAR}}
\newcommand{\theorist}{{\sc Theorist}}
\newcommand{\algep}{{\sc ALGEB}}
\newcommand{\algeb}{{\sc ALGEB}}
\newcommand{\amore}{{AMORE}}
\newcommand{\crep}{{CREP}}
\newcommand{\meataxe}{{\sc Meat-Axe}}
\newcommand{\moc}{{\sc MOC}}
\newcommand{\redux}{{\sc ReDuX}}
\newcommand{\sisyphos}{{\sc Sisyphos}}
\font\logo=cmss10 scaled \magstep4
\newdimen\r \r=3.14 pt \multiply\r by \magstep3\divide\r by 1000
\newbox\LiEbox \newbox\liepic
\setbox\LiEbox=\hbox{{\logo L\kern-.25em\raise\r\hbox{\i}\kern-0.03em E} }
\def\LiE.{\unhcopy\LiEbox}
\font\smlogo=cmss10 \r=2.718pt
\newbox\smallLiEbox \newbox\liepic
\setbox\smallLiEbox=\hbox{{\smlogo L\kern-.25em\raise\r\hbox{\i}\kern-0.03em E}}
\def\lie{\unhcopy\smallLiEbox}

\renewcommand{\baselinestretch}{0.9}

%
%
%
%
\newcommand{\Z} {{\rm {\mbox{\protect\makebox[.2em][l]{\sf Z}\sf Z}}}}
\newcommand{\N} {{\rm {\mbox{\protect\makebox[.15em][l]{I}N}}}}
\newcommand{\bi}{\bibitem}
\newcommand{\1}{{\bf{1}}}
\newcommand{\2}{{\bf{2}}}
\newcommand{\4}{{\bf{4}}}
\newcommand{\5}{{\bf{5}}}
\newcommand{\6}{{\bf{6}}}
\newcommand{\7}{{\bf{7}}}
\newcommand{\8}{{\bf{8}}}
\newcommand{\9}{{\bf{9}}}
\newcommand{\0}{{\bf{0}}}

\begin{center}
{\Large\bf Algebraische Darstellung transzendenter Funktionen}
\end{center}
\smallskip
\begin{center}
Wolfram Koepf
\end{center}
\smallskip
\smallskip
Ich m"ochte in diesem Bericht algorithmische Methoden vorstellen,
die im wesentlichen in diesem Jahrzehnt Einzug in die Computeralgebra
gefunden haben. Die haupts"achlichen Ideen gehen auf Stanley \cite{Sta} und
Zeilberger \cite{Zei1}--\cite{Zei4} zur"uck, 
vgl.\ die Beschreibung \cite{Strehl1}, und haben ihre Wurzeln teilweise bereits
im letzten Jahrhundert (siehe z.\ B.\ \cite{Beke1}--\cite{Beke2}), 
gerieten aber auf Grund der Komplexit"at der auftretenden
Algorithmen wieder in Vergessenheit.

Eine der Kernfragen kann hierbei so formuliert zu werden:
Worin liegt der wesentliche Unterschied zwischen
der Exponentialfunktion $f(x)=e^x$ und beispielsweise
der Funktion $g(x)=e^x+|x|/10^{1000}$, der dazu f"uhrt,
da{\ss} alle Mathematiker die
Exponentialfunktion $f$ als elementare Funktion auffassen und nicht $g$,
obwohl sich diese beiden Funktionen auf einem gro"sen Bereich der reellen Achse
numerisch kaum unterscheiden?
%

Oder ein Beispiel aus der diskreten Mathematik: Warum erh"alt die
Fakult"atsfunktion $a_n=n!$ den Vorzug gegen"uber
$b_n=n!+n/10^{1000}$ oder irgendeiner anderen diskreten Funktion?

Obwohl die Beispiele die bekanntesten stetigen bzw.\ diskreten
{\bf transzendenten} (nicht-\linebreak algebraischen) Funktionen betreffen, ist 
die Antwort auf diese Fragen rein {\bf algebraischer Natur}:
Die Exponentialfunktion $f$ ist n"amlich 
charakterisiert durch jede der folgenden algebraischen Eigenschaften:
\begin{enumerate}
\item
$f$ ist stetig, und f"ur alle $x,y$ gilt $f(x+y)=f(x)\cdot f(y)$;
\item
$f$ ist differenzierbar, und $f'(x)=f(x)$ sowie $f(0)=1$;
\item
$f\in C^{\infty}$, $f(x)=\sum\limits_{n=0}^{\infty}a_n\,x^n$ mit
$a_0=1$, und f"ur alle $n\geq 0$ gilt $(n+1)\,a_{n+1}=a_n$;
\end{enumerate}
und die Fakult"atsfunktion wird charakterisiert durch jede der folgenden
algebraischen Eigenschaften:
\begin{enumerate}
\item[4.]
$a_0=1$, und f"ur alle $n\geq 0$ gilt $a_{n+1}=(n+1)\,a_n$;
\item[5.]
die erzeugende Funktion $f(x)\!=\!\sum\limits_{n=0}^{\infty}a_n x^n$
erf"ullt die Differentialgleichung $x^2 f'(x)+$\linebreak $(x-1)f(x)+1\!=\!0$
mit der Anfangsbedingung $f(0)=1$. 
%
\end{enumerate}
W"ahrend mir keine Methode bekannt ist, Funktionalgleichungen wie
Eigenschaft (1.) in der Computeralgebra zur Darstellung transzendenter
Funktionen einzusetzen, will ich aufzeigen, in welcher Form
die anderen Eigenschaften dazu geeignet sind.

Man beachte, da"s die erzeugende Funktion der Fakult"at nur
am Ursprung konvergiert, also als formale Reihe anzusehen ist.
D.\ h.\ insbesondere, da"s man eine geschlossene Darstellung der erzeugenden
Funktion nicht angeben kann. (Diese ist zumindest nicht eindeutig: Jede
Funktion der Funktionenschar $-\frac{e^{-1/x}}{x}({\rm Ei}\:(-1/x)+c)$
(s.\ \cite{AS}, Kapitel 5)
stellt die erzeugende Funktion der Fakult"at f"ur $x>0$ asymptotisch dar.)
Darauf kommt es aber auch gar nicht an:
Statt mit der erzeugenden Funktion selbst arbeitet man ohnehin viel besser
mit ihrer Differentialgleichung, welche rein algebraisch ist!
Genau dasselbe trifft auf die Exponentialfunktion und die Fakult"at selbst zu:
Statt mit diesen transzendenten Objekten arbeite man mit den dazu
"aquivalenten algebraischen Differential- bzw.\ Rekursionsgleichungen
(vgl.\ \cite{Buch}).

Die gegebenen Eigenschaften sind {\bf Strukturaussagen} f"ur die
betreffenden Funktionen. Bei der geringsten "Anderung geht diese
Struktur verloren. Beispielsweise kann die Funktion $g(x)=e^x+|x|/10^{1000}$
durch keine den Eigenschaften (1.)--(3.) analoge Vorschrift 
charakterisiert werden. Dagegen ist die Funktion $h(x)=e^x+x/10^{1000}$ 
beispielsweise durch die Differentialgleichung 
$(x-1)\,h''(x)-x\,h'(x)+h(x)=0$ mit den Anfangsbedingungen
$h(0)=1$ und $h'(0)=1+10^{-1000}$ gegeben.

Das Besondere (und Gemeinsame) an Exponential- bzw.\ Fakult"atsfunktion
besteht also darin, da"s diese eine homogene lineare Differentialgleichung
und jene eine homogene lineare Rekursionsgleichung erf"ullt, wobei
Differentialgleichung sowie Rekursionsgleichung Polynomkoeffizienten
haben und beide erster Ordnung sind. 

Verallgemeinern wir diesen Sachverhalt nun zun"achst auf stetige Funktionen
einer Variablen und nennen eine Funktion $f(x)$ {\bf holonom}, falls sie
eine homogene lineare Differentialgleichung mit Polynomkoeffizienten 
in $x$ erf"ullt.
Stanley \cite{Sta} zeigte, da"s Summe und Produkt holonomer Funktionen
sowie die Komposition mit algebraischen Funktionen wieder holonome
Funktionen liefern. Beke \cite{Beke1}--\cite{Beke2} hat bereits vor 
100 Jahren Algorithmen beschrieben, mit welchen die Differentialgleichung
f"ur Summe bzw.\ Produkt von $f$ und $g$ aus den Differentialgleichungen f"ur
$f$ und $g$ bestimmt werden k"onnen!

Analog nennt man eine diskrete Funktion
$a_n$ holonom, falls sie eine homogene lineare Rekursionsgleichung 
mit Polynomkoeffizienten in $n$ erf"ullt. Summe und Produkt holonomer 
diskreter Funktionen sind wieder holonom, und es gibt Algorithmen zur 
Berechnung der entsprechenden Rekursionen (s.\ \cite{SZ}, \cite{KS}).

Was haben wir hiermit nun gewonnen? Ignorieren wir einmal, da"s
$e^x,$ $\sin x, \cos x, \arctan x,$ $\arcsin x$ etc.\ transzendente Funktionen 
sind, und stellen wir lediglich ihre holonomen Differentialgleichungen
$f'=f$, $f''=-f$, $f''=-f$, $(1+x^2)f''+2xf'=0$, $(x^2-1)f''+xf'=0$ etc.\
in Rechnung, so k"onnen wir nun aus diesen Differentialgleichungen
mit reiner Polynomarithmetik (wir brauchen eigentlich nur lineare Algebra,
vgl.\ \cite{SZ}, \cite{KS})
holonome Differentialgleichungen f"ur Summen und Produkte solcher Funktionen,
beispielsweise f"ur $f(x)=\arcsin^2 x$ (n"amlich $(x^2-1)f'''+3xf''+f'=0$),
erzeugen. Doch damit nicht genug: F"ur die Koeffizienten $a_n$ der formalen 
Potenzreihe von $\arcsin^2 x=\sum\limits_{n=0}^\infty a_n x^n$ folgt dann
automatisch die holonome Rekursionsgleichung $n(1 + n) (2 + n)a_{n+2}=n^3a_n$,
welche (gl"ucklicherweise) nur die zwei Terme $a_{n+2}$ und $a_n$
enth"alt, sich daher l"osen l"a"st und zu der Darstellung 
\[
\arcsin^2 x=\sum_{n=0}^\infty
\frac{4^n\,n!^2}{(1 + n)\,(1 + 2n)!}x^{2n+2}
\]
f"uhrt (vgl.\ \cite{GKP}, \cite{Wil1}, \cite{Koe92}--\cite{Koe93}).

Oder notieren wir uns lediglich die holonomen Rekursionen
rationaler Funktionen sowie die der Funktionen
\begin{equation}
(mn+b)!\;,
\quad\quad
\frac{1}{(mn+b)!}\;
\quad\quad(m\in\Z)
\quad\quad
\mbox{sowie}
\quad\quad
a^n
\label{eq:zulaessige diskrete Funktionen}
\end{equation}
(d.\ h.\ wir geben sie unserem favorisierten Computeralgebrasystem 
in geeigneter Form bekannt), so lassen sich nun f"ur alle m"oglichen durch 
Addition und Multiplikation erzeugten Funktionen holonome Rekursionen herleiten,
beispielsweise die beiden Rekursionen
\begin{equation}
(n-k+1)^2 F(n+1,k)-(1+n)^2 F(n,k)=0
\label{eq:combnk^2a}
\end{equation}
und
\begin{equation}
(k+1)^2 F(n,k+1)-(n-k)^2 F(n,k)=0
\label{eq:combnk^2b}
\end{equation}
f"ur $F(n,k)={{n}\choose{k}}^2$. Diese k"onnen zwar auch leicht direkt
aus der Darstellung von $F(n,k)$ abgelesen werden,
aber das vorgestellte Verfahren kann problemlos auf die kompliziertere
Funktion $F(n,k)=\frac{n!+k!^2}{k}$ angewandt werden und 
f"uhrt dann auf die Rekursionen
\[
0=n F(n+2,k) -(1+3n+n^2) F(n+1,k) +(1+n)^2 F(n,k)
\]
und
\[
0=k (2 + k)^2 F(n,k+2) -(1 + k) (1 + 3 k + k^2) (3 + 3 k + k^2) F(n,k+1)
+ k (1 + k)^3 F(n,k)
\;.
\]
%
%
%
%
%
%
%
%
%
Eine wichtige Fragestellung der Kombinatorik ist,
zu einer gegebenen Funktion $F(n,k)$, 
ausgedr"uckt als Produkt von Termen der Form 
(\ref{eq:zulaessige diskrete Funktionen}), die unendliche Summe
\[
s(n)=\sum_{k} F(n,k)
\]
zu berechnen, wobei "uber alle $k\in\Z$ zu summieren ist. 
(In der Praxis sind jedoch h"aufig nur endlich viele Terme von Null 
verschieden.) Ist nun $F(n,k)$ ein {\bf hypergeometrischer Term},
\mbox{d.\ h.,} sind sowohl $F(n+1,k)/F(n,k)$ als auch $F(n,k+1)/F(n,k)$
rational bzgl.\ beider Variablen $n$ und $k$ -- dies trifft beispielsweise
wegen (\ref{eq:combnk^2a})--(\ref{eq:combnk^2b})
f"ur $F(n,k)={{n}\choose{k}}^2$ zu --
so findet der (schnelle) {\bf Zeilberger-Algorithmus} (\cite{Zei2}, s.\ auch
\cite{Koornwinder} und \cite{PS}) eine holonome
Darstellung, d.\ h.\ eine holonome Rekursion, f"ur $s(n)$.
Dieser Algorithmus baut auf dem von Gosper \cite{Gos}
gefundenen Entscheidungsalgorithmus f"ur die unbestimmte Summation auf.
In unserem Beispielfall findet der Zeilberger-Algorithmus f"ur
$s(n)=\sum\limits_{k} {{n}\choose{k}}^2=\sum\limits_{k=0}^n {{n}\choose{k}}^2$ 
die holonome Rekursion 
\[
(1+n)\,s(n+1)=2(1+2n)\,s(n)\;,
\]
welche (gl"ucklicherweise) wieder nur zwei Terme hat. Daher bekommen
wir die Darstellung
\[
s(n)=\sum\limits_{k=0}^n {{n}\choose{k}}^2=\frac{(2n)!}{n!^2}
\;.
\]
Im allgemeinen wird die resultierende Rekursion bzw.\ Differentialgleichung
nat"urlich mehr als zwei
Terme enthalten. Aber auch dann enth"alt diese zum einen eine interessante
Strukturinformation (z.\ B.\ "uber die Orthogonalit"at eines Polynomsystems)
und kann zudem eine f"ur numerische Zwecke n"utzliche Vorschrift darstellen 
(vgl.\ \cite{Deufl1}--\cite{Deufl2}).

Die gefundene Strukturinformation kann insbesondere zur {\bf Identifikation} 
transzendenter Funktionen herangezogen werden.
Um z.\ B.\ die Identit"at 
\[
\sum_{k=0}^n
{{n}\choose{k}}^3
=
\sum_{k=0}^n
{{n}\choose{k}}^2 {{2k}\choose{n}}
\]
(vgl.\ \cite{Strehl2}) zu "uberpr"ufen -- diese ist nicht trivial, da die
beiden Summanden {\bf nicht} dieselben Rekursionen bzgl.\ $k$ und $n$
erf"ullen! -- brauchen wir nur zu zeigen, da"s beide Summen derselben Rekursion
\[
(n+2)^2 s(n+2)-(16 + 21 n + 7 n^2)s(n+1)-(n+1)^2 s(n)=0
\]
gen"ugen -- dies macht der Zeilberger-Algorithmus -- sowie dieselben 
Anfangswerte $s(0)=1$ und $s(1)=2$ haben (dies ist trivial).

Als weiteres Beispiel betrachten wir die Funktion 
($\alpha,\beta,\gamma\in\N_0$)
\begin{eqnarray*}
V(\alpha,\beta,\gamma)
&=&
(-1)^{\alpha+\beta+\gamma}\cdot
\frac{\Gamma(\alpha\!+\!\beta\!+\!\gamma\!-\!d/2)
\Gamma(d/2\!-\!\gamma)\Gamma(\alpha\!+\!\gamma\!-\!d/2)
\Gamma(\beta\!+\!\gamma\!-\!d/2)}
{\Gamma(\alpha)\Gamma(\beta)\Gamma(d/2)\Gamma(\alpha+\beta+2\gamma-d)
(m^2)^{\alpha+\beta+\gamma-d}}
\\\\&&\cdot
\; _2 F_1\left.
\!\!
\left(
\!\!\!\!
\begin{array}{c}
\multicolumn{1}{c}{\begin{array}{cc} \alpha+\beta+\gamma-d
\;, & \alpha+\gamma-d/2 \end{array}}\\[1mm]
\multicolumn{1}{c}{ \alpha+\beta+2\gamma-d}
            \end{array}
\!\!\!\!
\right| z\right)
\end{eqnarray*}
($ _2 F_1$ stellt hierbei die Gau"ssche hypergeometrische Reihe dar, s.\
z.\ B.\ \cite{AS}, Kapitel 15),
welche bei der Berechnung von Feynman-Diagrammen eine Rolle spielt 
\cite{FT1}, f"ur die wir die Rekursion
\begin{eqnarray*}
0&=&
  \left( 2\,\alpha \! -\! d\! +\! 2\,\gamma  \right) \,
     \left( 2\,\alpha \! +\! 2\,\beta\!  -\! d \!+\! 2\,\gamma  \right) \,
     \left( 2\! +\! 2\,\alpha \! +\! 2\,\beta\!  - \!d \!+\! 2\,\gamma  \right) \,
     { V}(\alpha ,\beta ,\gamma )
\\&&
-
    2\alpha \left( 2\! + \!2\alpha \! +\! 2\beta \! - \!d\! +\! 2\gamma  \right)
     {m^2}
\\&&\hspace*{1cm}
\cdot
\left( -2\alpha \! - \!2\beta \! + \!2d \!- 4\gamma \! +\! 2z \!+
       \!4\alpha z \!+\! 2\beta z\! -\! 3dz \!+ \!4\gamma z \right)
     { V}(1 \!+ \!\alpha ,\beta ,\gamma )
\\&&+
    8\,\alpha \,\left( 1 + \alpha  \right) \,
     \left( 1 + \alpha  + \beta  - d + \gamma  \right) \,{m^4}
     \left( z-1 \right) \,z\,{ V}(2 + \alpha ,\beta ,\gamma )
\end{eqnarray*}
sowie analoge Rekursionen bzgl.\ der Variablen $\beta$ und $\gamma$
erhalten. Diese k"onnen dann zur numerischen Berechnung herangezogen werden.

Zeilberger 
betrachtete in \cite{Zei1} die allgemeinere Situation von Funktionen $F$
mehrerer diskreter und stetiger Variabler. Sind es $d$ Variablen
und hat man $d$ (im wesentlichen unabh"angige) m"oglicherweise gemischte
homogene partielle Differential-Differenzengleichungen mit 
Polynomkoeffizienten (bzgl.\ aller Variablen) f"ur $F$, nennen wir $F$
ein holonomes System (vgl.\ \cite{Bernstein1}--\cite{Bj}). 
Dann legen diese Gleichungen $F$ zusammen mit
geeigneten Anfangswerten bereits eindeutig fest.

Insbesondere gilt dies also, wenn das gegebene System holonomer
Gleichungen separiert ist, d.\ h., wenn in jeder der Gleichungen
nur Ableitungen bzgl.\ einer der stetigen Variablen bzw.\ nur Shifts
bzgl.\ einer der diskreten Variablen vorkommen.
Beispielsweise bilden die Legendre-Polynome $F(n,x)=P_n(x)$
(\cite{AS}, Kapitel 22) auf Grund ihrer Differentialgleichung
\begin{equation}
(x^2-1)F''(n,x)+2xF'(n,x)-n(1+n)F(n,x)=0
\label{eq:Legendre1}
\end{equation}
sowie ihrer Rekursionsgleichung
\begin{equation}
(n+2)F(n+2,x)-(3 + 2 n) x  F(n+1,x)+(n+1)F(n,x)=0
\label{eq:Legendre2}
\end{equation}
zusammen mit den Anfangsbedingungen
\[
F(0,0)=1
\;,\quad
F(1,0)=0
\;,\quad
F'(0,0)=0
\;,\quad
F'(1,0)=1
\]
%
%
%
ein holonomes System. Gleichungen (\ref{eq:Legendre1}) und (\ref{eq:Legendre2})
stellen also eine hinreichende algebraische, ja sogar polynomiale,
Struktur zur Beschreibung von $P_n(x)$ dar.

Fa"st man die auftretenden (partiellen) Differentiationen und 
Indexverschiebungen als Operatoren und die 
Differenzen-Differentialgleichungen als Operatorengleichungen auf,
so stellen diese ein {\bf polynomiales Gleichungssystem} in
einem nichtkommutativen Polynomring dar. Ist n"amlich $x$ eine
stetige Variable und $D_x$ der zugeh"orige Ableitungsoperator,
so ist wegen der Produktregel $D_x (x f)-x D_x f=f$, und
folglich hat man den Kommutator $D_x x-xD_x=1$.
Ist andererseits $k$ eine diskrete Variable und $K$ der zugeh"orige
Indexverschiebungsoperator $K a_k=a_{k+1}$ (Aufw"artsshift),
so gilt $K (k a_k)-k K a_k=$ \mbox{$(k+1) a_{k+1}-k a_{k+1}=a_{k+1}=Ka_k$,}
und folglich gilt die Kommutatorregel $K k-k K=K$.
Entsprechendes gilt f"ur die restlichen Variablen, w"ahrend
alle anderen Kommutatoren verschwinden.

Das Umformen eines durch gemischte Differenzen-Differentialgleichungen 
gegebenen holonomen Systems stellt sich
in dem betrachteten nichtkommutativen Polynomring als ein
polynomiales Eliminationsproblem dar,
welches mit nichtkommutativen Gr"obnerbasen-Methoden gel"ost werden
kann (\cite{BW}, \cite{Gal}, \cite{Weispfenning}, \cite{Zei1}, \cite{Zei4},
\cite{Ta1}--\cite{Ta3}). 
Als Beispiel diene $F(n,k)={{n}\choose{k}}$.
Hierf"ur gilt die Pascalsche Dreiecksbeziehung 
$F(n+1,k+1)=F(n,k)+F(n,k+1)$ sowie die reine Rekursion
$(n+1-k) F(n+1,k)-(n+1) F(n,k)=0$ bzgl.\ $n$, welche in Operatornotation
$(KN-1-K) F(n,k)=0$ sowie $((n+1-k)N-$ $(n+1))F(n,k)=0$ lauten, wobei
$N F(n,k)=F(n+1,k)$ den Verschiebungsoperator bzgl.\ $n$ bezeichne.
Somit erh"alt man das Polynomsystem
\[
KN-1-K=0
\quad\quad\mbox{sowie}\quad\quad
(n+1-k)N-(n+1)=0
\;.
\]
Die Gr"obnerbasis des erzeugten Linksideals bzgl.\ der 
Termordnung $(k,n,K,N)$ (lexikographisch) ergibt sich mit \cite{AM} zu
%
%
%
\[
\Big\{(k+1) K +k-n,(n+1-k)N-(n+1),KN-1-K\Big\}
\;,
\]
d.\ h.\ also, da"s auf diesem Wege {\bf automatisch} die reine Rekursion 
\[
(k+1)\,F(n,k+1) +(k-n)\,F(n,k)=0
\]
bzgl.\ $k$ erzeugt wurde.

Als weiteres Beispiel betrachte ich die Legendre-Polynome, f"ur die die
Beziehungen (\ref{eq:Legendre1})--(\ref{eq:Legendre2}) gelten. Hier haben
wir also das Polynomsystem ($D:=D_x$)
\[
(x^2-1)D^2+2xD -n(1+n)=0
\quad\mbox{sowie}\quad
(n+2)N^2-(3+2n)x N+(n+1)=0
\;.
\]
Die Gr"obnerbasis des erzeugten Linksideals bzgl.\ der
Termordnung  $(D,N,n,x)$ ergibt sich zu 
\[
\Big\{
(x^2-1)D^2+2xD -n(1+n),
\]
\[
(1+n)ND-(1+n)xD-(1+n)^2,
\]
\begin{equation}
(x^2-1)ND-(1+n)xN+(1+n),
\label{eq:Legendre3}
\end{equation}
\begin{equation}
(1+n)(x^2-1)D-(1+n)^2 N+x(1+n)^2
\label{eq:Legendre4}
\end{equation}
\[
(n+2)N^2-(3+2n)x N+(n+1)
\Big\}
\;,
\]

\pagebreak\noindent
wobei ich der besseren Lesbarkeit halber die Operatoren $D$ und $N$ wieder
rechts positioniert habe.
Hier wurden also $D$-Potenzen weitestgehend eliminiert, und Gleichungen
(\ref{eq:Legendre3})--(\ref{eq:Legendre4}) entsprechen den Beziehungen
\begin{equation}
(x^2-1)P_{n+1}'(x)=(1+n)\,(x P_{n+1}(x)-P_n(x)) 
\label{eq:Legendre5}
\end{equation}
\begin{equation}
(x^2-1)P_{n}'(x)=(1+n)\,(P_{n+1}(x)-x P_n(x))
\label{eq:Legendre6}
\end{equation}
zwischen den Legendre-Polynomen und ihren ersten Ableitungen.
Man sieht also, da"s auf diesem Wege neue Beziehungen (zwischen den
Binomialkoeffizienten bzw.\ zwischen den Ableitungen der Legendre-Polynome)
{\bf hergeleitet} wurden.

Analog lassen sich mit dieser Methode Rekursionen f"ur {\bf holonome Summen}
herleiten. Betrachten wir beispielsweise 
\[
s(n)=\sum_{k=0}^n F(n,k)=\sum_{k=0}^n {{n}\choose{k}} P_n(x)
\;,
\]
dann findet man mit dem Produktalgorithmus zun"achst die holonomen Rekursionen
\[
(n-k+1) F(n+1,k)-(1+n) F(n,k) = 0
\]
sowie
\[
(2+k)^2 F(n,k+2)-(3+2k)(n-k-1)x F(n,k+1)+(n-k)(n-k-1) F(n,k)=0
\]
f"ur den Summanden $F(n,k)$. In der Gr"obnerbasis des von den zugeh"origen
Polynomen
\[
(n-k+1) N-(1+n)
\]
sowie
\[
(2+k)^2 K^2-(3+2k)(n-k-1)x  K +(n-k)(n-k-1) 
\]
bzgl.\ der Termordnung $(k,n,K,N)$ erzeugten Linksideals liegt das 
$k$-freie Polynom
\[
(2+n)^2 K^2 N^2 - K (2 + n) (3 + 2 n) (K + x) N + 
(1 + n) (2 + n) (1 + K^2  + 2 K x)
\;,
\]
welches einer $k$-freien Rekursion f"ur $F(n,k)$ entspricht.
Da bei der Summation "uber $k\in\Z$ die verschobenen Summen
\[
s(n)=\sum_{k} F(n,k)=\sum_{k} F(n,k+1)=\sum_{k} F(n,k+2)
\]
alle dieselbe Summenfunktion $s(n)$ liefern,
liefert die Substitution $K:=1$ die g"ultige holonome Rekursion 
\[
(2+n) s(n+2) - (3 + 2 n) (1 + x) s(n+1) + 2 (1 + n) (1 + x) s(n)=0
\]
f"ur $s(n)$.

Die gezeigte Methode ist zur Generierung von Identit"aten im Prinzip
universell einsetzbar, ben"otigt jedoch
den komplizierten Apparat des (nichtkommutativen) Buchberger-Algorithmus 
und erbt die damit verbundenen Nachteile. Eine wesentliche H"urde stellt
die Komplexit"at bei Problemen mit vielen Variablen dar.

Interessiert man sich wie bei obigem Beispiel 
(\ref{eq:Legendre5})--(\ref{eq:Legendre6})
f"ur die Erzeugung von Identit"aten zwischen
den Ableitungen $F^{(j)}(n+k,x)\;(j,k\in\N_0)$ 
eines holonomen Systems $F(n,x)$, 
mu"s dies aber nicht unbedingt sein.
Es geht in vielen F"allen bereits mit {\bf linearer} Algebra! Dazu mu"s
man aber mehr Information hineinstecken.
Die geeignete "uber die holonomen Beziehungen hinausgehende Information 
besteht in einer Beziehung der Form
\[
F'(n,x)=\sum_{k=0}^m r_k(n,x) F(n+k,x)
\;,
\]
mit rationalen Funktionen $r_k$ bzgl.\ $n$ und $x$,
einer {\bf Ableitungsregel} f"ur $F(n,x)$ also, sofern erh"altlich. 
Es zeigt sich, da"s in der Praxis holonome Systeme (wie beispielsweise
Systeme orthogonaler Polynome etc., s.\ z.\ B.\ \cite{AS}, 22.8, 
und \cite{Koeortho}) so strukturstark sind, 
da"s eine Ableitungsregel verf"ugbar ist. Man kann zeigen, da"s
Summe und Produkt solcher Systeme in der Regel auch wieder holonome
Systeme mit Ableitungsregel darstellen \cite{Koeortho},
und Abh"angigkeiten zwischen den Ableitungen $F^{(j)}(n+k,x)\;(j,k\in\N_0)$
k"onnen mit reiner linearer Algebra gefunden werden.

Auf diese Weise wurde z.\ B.\ die Beziehung
\begin{eqnarray*}
P_n^{(\alpha,\beta)}(x)
&=&
-\frac{2 (\alpha + n) (\beta + n)}
{(\alpha + \beta + n) (\alpha + \beta + 2 n) (\alpha + \beta + 2 n+1)}\,
{P_{n-1}^{(\alpha,\beta)}}'(x)
\\&&
+\frac{2 (\alpha - \beta)}
{(\alpha + \beta + 2 n) (\alpha + \beta + 2 n+2)}\,
{P_{n}^{(\alpha,\beta)}}'(x)
\\&&
+
\frac{2 (\alpha + \beta + n+1)}
{(\alpha + \beta + 2 n+1) (\alpha + \beta + 2 n+2)}\,
{P_{n+1}^{(\alpha,\beta)}}'(x)
\end{eqnarray*}
f"ur die Jacobi-Polynome $P_n^{(\alpha,\beta)}(x)$ (\cite{AS}, Kapitel 22)
automatisch erzeugt.
Hierbei war die Zielsetzung, da"s die auftretenden Koeffizientenfunktionen
von ${P_{n+k}^{(\alpha,\beta)}}'(x)$
nicht von $x$ abh"angen sollen. Dies ist f"ur Fragestellungen aus der
Spektralapproximation (s.\ \cite{CHQZ}, \S~2.3.2) von Bedeutung.

Zuletzt will ich darauf verweisen, da"s die vorliegende Beschreibung
nat"ur\-lich keinen Anspruch auf Vollst"andigkeit erheben kann.
Ich konnte weder auf den Gosper-Algorithmus \cite{Gos}
noch auf die Wilf-Zeilberger-Theorie der WZ-Paare und rationalen Zertifikate
eingehen \cite{Wilf}--\cite{WZ4}.
Auch Petkov\u seks Algorithmus
\cite{P}, der alle hypergeometrischen Terml"osungen holonomer
Rekursionsgleichungen berechnet, konnte keine Ber"ucksichtigung finden.

Die in diesem Artikel durchgef"uhrten Rechnungen wurden mit \mathematica\
und \reduce\ durchgef"uhrt. Mein Dank gilt Prof.\ Peter Deu"flhard, der mich
ermutigt hat, mich mit dem vorliegenden Thema zu besch"aftigen,
Herbert Melenk, mit dem ich wichtige Gespr"ache "uber
nichtkommutative Gr"obnerbasen f"uhren konnte, sowie Jochen Fr"ohlich,
der mich auf die Spektralapproximation aufmerksam machte.

\end{document}